\documentstyle[12pt]{article}
\def\beq{\begin{equation}}
\def\eeq{\end{equation}}
\def\rc{\check{R}}
\newcommand{\mdd}[4]
{\left(\begin{array}{rr} #1 & #2 \\ #3 & #4 \end{array}\right)}
\newcommand{\stt}{$(**)$}
\newtheorem{theo}{Theorem}
\newtheorem{cor}{Corollary}
\newtheorem{lem}{Lemma}
\def\IR{\relax{\rm I\kern-.18em R}}
\font\cmss=cmss10 \font\cmsss=cmss10 at 7pt
\def\IZ{\relax\ifmmode\mathchoice
{\hbox{\cmss Z\kern-.4em Z}}{\hbox{\cmss Z\kern-.4em Z}}
{\lower.9pt\hbox{\cmsss Z\kern-.4em Z}}
{\lower1.2pt\hbox{\cmsss Z\kern-.4em Z}}\else{\cmss Z\kern-.4em Z}\fi}
\def\inbar{\,\vrule height1.5ex width.4pt depth0pt}
\def\IC{\relax\hbox{$\inbar\kern-.3em{\rm C}$}}
\def\ID{\relax{\rm I\kern-.18em D}}
\def\IE{\relax{\rm I\kern-.18em E}}
\def\IF{\relax{\rm I\kern-.18em F}}
\def\IG{\relax\hbox{$\inbar\kern-.3em{\rm G}$}}
\def\IH{\relax{\rm I\kern-.18em H}}
\def\II{\relax{\rm I\kern-.18em I}}
\def\IK{\relax{\rm I\kern-.18em K}}
\def\IL{\relax{\rm I\kern-.18em L}}
\def\IM{\relax{\rm I\kern-.18em M}}
\def\IN{\relax{\rm I\kern-.18em N}}
\def\IO{\relax\hbox{$\inbar\kern-.3em{\rm O}$}}
\def\IP{\relax{\rm I\kern-.18em P}}
\def\IQ{\relax\hbox{$\inbar\kern-.3em{\rm Q}$}}
\def\IGa{\relax\hbox{${\rm I}\kern-.18em\Gamma$}}
\def\IPi{\relax\hbox{${\rm I}\kern-.18em\Pi$}}
\def\ITh{\relax\hbox{$\inbar\kern-.3em\Theta$}}
\def\IOm{\relax\hbox{$\inbar\kern-3.00pt\Omega$}}

\def\Z{\IZ}

\def\C{\IC}
\begin{document}
\begin{titlepage}
\begin{center}
{}~\hfill ETH-TH/95-11\\
{}~\hfill March 1995
\end{center}
\vskip 1.5in
\begin{center}
{\large \bf Representations of knot groups
and Vassiliev invariants}
\vskip 0.5in
Daniel Altschuler\\
{\em Institut f\"ur Theoretische Physik\\
ETH-H\"onggerberg\\
CH-8093 Z\"urich, Switzerland\\
altsch@itp.phys.ethz.ch}
\end{center}
\vskip .5in
\begin{abstract}
We show that the number of homomorphisms
from a knot group to a finite group $G$
cannot be a Vassiliev invariant, unless it is constant
on the set of $(2,2p+1)$ torus knots. In several cases,
such as when $G$ is a dihedral or symmetric group, this
implies that the number of homomorphisms is not a
Vassiliev invariant.
\end{abstract}
\end{titlepage}
Recently, the Vassiliev knot invariants \cite{vassil}, also known
as invariants of finite type, have attracted a lot of interest.
The main properties of these invariants can be found in \cite{intro}.
In some sense the Vassiliev invariants include all the invariants
associated to quantum deformations of Lie algebras.
However, many classical knot invariants, such as
the crossing number, the braid index and the signature are
not Vassiliev invariants \cite{bl,dean,trapp}.
Here we show that another
well-known invariant, the number of homomorphisms
(representations) of a knot group in a finite group $G$,
is not a Vassiliev invariant, for several classes of groups $G$.

Let $K$ be a knot, and denote by $G(K)=\pi_1(S^3-K)$ the
knot group. Let $G$ be a finite group, and
$C$ a subset of $G$ stable under
conjugation: $gCg^{-1}=C$ for all $g\in G$.
Let $[K,G,C]$ be the number of homomorphisms
$G(K)\rightarrow G$ such that the image of a meridian of $K$
is contained in $C$.
In this note, we prove that if $K\mapsto [K,G,C]$ is not constant
on the set of $(2,2p+1)$ torus knots, it is
not an invariant of finite type, see Corollary 1 below.
For a given pair $(G,C)$, it is a fairly simple
matter to check whether this condition is satisfied. We illustrate
this with various examples, including the dihedral and
symmetric groups, for which we are able to conclude that
$[K,G,C]$ is not an invariant of finite type.

The dihedral representations of knot groups have been first
considered by Fox \cite{fox}. A states model definition of
$[K,G,C]$ is given in \cite{hj}.

One preliminary remark is in order.
Let $G'$ be the commutator subgroup of $G$. If
$\varphi : G(K)\rightarrow G$ is surjective, then $G/G'$ is
cyclic since $G(K)/G(K)' = \Z$. Therefore if $G/G'$ is not
cyclic, $[K,G,C]=[K,H,C\cap H]$ for some proper subgroup $H$
of $G$. Thus we can always assume that $G/G'$ is cyclic.

Recall the definition of an invariant of finite type.
Consider a set $X=\{c_1,\ldots,c_N\}$ of
$N$ crossings of a regular projection of a knot $K$. For
each map $\epsilon : X \rightarrow \{\pm 1\}$, let $K_\epsilon$
be the knot obtained by switching the crossing $c_i$
if $\epsilon(c_i)=-1$, for all $i=1,\ldots,N$.
Let $v(K)$ be a $\C$-valued invariant of knots.
Then $v(K)$ is said to be
of type $n$, where $n$ is a non-negative integer, if for all $K$ and
all sets $X$ of $N$ crossings with $N>n$,
\beq
\sum_\epsilon \prod_{i=1}^{N} \epsilon(c_i) v(K_\epsilon) = 0,
\eeq
where summation is over all maps $\epsilon : X\rightarrow\{\pm1\}$.
An invariant is of finite type if there exists $n$ such that
it is an invariant of type $n$.

Our starting point is a bound
for the positive number $[K,G,C]$.
It is a particular case of a more general theorem of
Turaev, see \cite{turaev}, p. 114. Using Turaev's
definitions, one could derive this bound by
proving that the category of representations of the
quantum double of $G$ is a unitary modular category.
However we will give a direct, elementary proof.
Let $c$ be the number of elements in $C$.
\begin{theo}
If $K$ may be obtained as the closure of a braid on $k$ strands,
then $[K,G,C]\leq c^k$.
\end{theo}
{\em Proof.} Let $V$ be the complex vector space with basis
$C$. We construct a representation of the
braid group $B_k$ on $V^{\otimes k}$. Define an
automorphism $\rc$ of $V^{\otimes 2}$ by
\beq
\rc(a\otimes b) = aba^{-1} \otimes a,
\label{w1}
\eeq
where $a,b\in C$.
Its inverse is given by
\beq
\rc^{-1}(a\otimes b) = b \otimes b^{-1}ab.
\label{w2}
\eeq
The desired representation is
$\rho : b_i \mapsto \rc_{i,i+1}$,
where $b_i$, $i=1,\ldots,k-1$ are the generators of $B_k$ and
\beq
\rc_{i,i+1}=
1\otimes 1\otimes\cdots\otimes\rc\otimes\cdots\otimes 1,
\eeq
where $\rc$ acts in the $i$\/-th and $i+1$\/-th copies of $V$.
By comparing (\ref{w1}) and (\ref{w2}) with the
relations in the Wirtinger presentation of $G(K)$, it is
easy to see that $[K,G,C]={\rm Tr} \,\rho(\beta)$, where
$\beta\in B_k$ is a braid such that $K$ is equivalent to
the closure of $\beta$.

Now observe that in the basis $C\otimes C$, $\rc$
is a permutation matrix. Thus for any $\beta\in B_k$,
$\rho(\beta)$ is also a permutation matrix.
Since the trace of a permutation matrix is always less or equal
than the trace of the identity matrix of the same size,
this proves the theorem. $\Box$

Remark that if all elements of $C$ commute among themselves,
then $\rc=\rc^{-1}$ and $[K,G,C]$ is a constant.

For any integer $p$, let $K_p$ denote the $(2,2p+1)$ torus knot.
Dean \cite{dean} and Trapp \cite{trapp}
independently proved the following result.
\begin{theo}
If $v(K)$ is an invariant of type $n$, then $v(K_p)$, as a
function of $p$, is a polynomial of degree at most $n$.
\end{theo}
\begin{cor}
If the function $p\mapsto [K_p,G,C]$ is not constant,
$[K,G,C]$ is not an invariant of finite type.
\end{cor}
{\em Proof.}
Applying theorem 1 to the case $K=K_p$, we find
$[K_p,G,C]\leq c^2$, because every $K_p$ may be realised
as the closure of a braid on two strands. By theorem 2,
if $[K_p,G,C]$ is of finite type, it is a polynomial in $p$,
so being bounded, it must be a constant. $\Box$

\begin{lem}
The following conditions are equivalent:
\begin{itemize}
\item[$(*)$] $[K_p,G,C]$ is not constant.
\item[\stt] There exist
$a,b\in C$, $a\neq b$, such that $(ab)^p a = b(ab)^p$ for
some integer $p$.
\end{itemize}
\end{lem}
{\em Proof.}
An easy induction shows that $G(K_p)$ has a presentation
with two generators $g,h$ and one relation $(gh)^p g = h(gh)^p$.
Therefore, $[K_p,G,C]$ is the number of elements in
\beq
{\cal C}_p = \{ (a,b)\in C\times C \, | \, (ab)^p a = b(ab)^p\}.
\eeq
For all $p$,
${\cal C}_p\supset \{ (a,a) \, | \, a\in C\} = {\cal C}_0$.
$\Box$

It is clear that if the pair $(G,C)$ satisfies these conditions,
then so does any pair $(G_1,C_1)$ with $G_1\supset G$, $C_1\supset C$.

Let $n>1$ be an integer. The dihedral group $D_n$ of order $2n$
is a group with two generators $r_1,r_2$ and the relations
$r_1^2 = r_2^2 = (r_1r_2)^n = 1$. Note that $D_n/D'_n$ is cyclic
if and only if $n$ is odd. More generally, a finite
Coxeter group $W$ of rank $\ell$ is a finite group with generators
$r_1,r_2,\ldots,r_\ell$ and relations
$(r_i r_j)^{m_{ij}}=1$, where $m_{ij}=m_{ji}$,
$m_{ij}\geq 2$ if $i\neq j$, and $m_{ii}=1$ for all $i$.
Note that any Coxeter group $W$
of rank $\ell\geq 2$ contains a dihedral group $D_n$.
Finite Coxeter group have been classified \cite{bou}. They
are Weyl groups of finite-dimensional semisimple Lie algebras,
dihedral groups, or two other groups $H_3$ and $H_4$.
The symmetric group $S_\ell$ is realized as the Weyl group of the
Lie algebra $sl_\ell$, the generators $r_i$ are the transpositions
$(i,i+1)$. The two transpositions $(i,i+1)$ and $(i+1,i+2)$ generate
a subgroup $D_3$.

\begin{lem}
If $n$ is odd and $C_n$ is the set of conjugates of the generators
$r_1,r_2$, $[K,D_n,C_n]$ is not
an invariant of finite type.
\end{lem}
{\em Proof.}
Take $p=(n-1)/2$. Then it is easy to check that
$a=r_1$, $b=r_2$, satisfy the condition \stt~above. $\Box$

Apart from the groups containing dihedral groups,
there is another class of groups for which we can easily show that
$[K,G,C]$ is not an invariant of finite type. They are the
groups containing $SL(2,\Z_m)$ or $PSL(2,\Z_m)$. Consider
$A,B\in SL(2,\Z_m)$ given by
\beq
A=\mdd{1}{1}{0}{1}, \;\; B=\mdd{1}{0}{-1}{1}.
\label{defab}
\eeq
It is again easy to check that $A$ and $B$ satisfy \stt~with
$p=1$. The reason for this is that the group of the trefoil has a
representation in $SL(2,\Z)$, since it has a presentation
with generators $x,y$ and one relation $x^2=y^3$.
Thus we get
\begin{lem}
If $G=SL(2,\Z_m)$ and $C$ is the set of conjugates of $A$ and $B$,
then $[K,G,C]$ is not an invariant of finite type.
\end{lem}
{}From these observations we deduce the next results on
some important families of groups.
\begin{theo}
In the following cases $[K,G,C]$ is not
an invariant of finite type:
\begin{itemize}
\item[{\rm (a)}]
$G$ is a finite Coxeter group which contains
a subgroup
$D_n$ with odd $n$, and $C$ contains the conjugates
of the generators $r_1,r_2$ of $D_n$.
\item[{\rm (b)}]
$G$ is the symmetric group $S_n$, $n\geq 3$, and
$C$ is the set of transpositions.
\item[{\rm (c)}]
$G$ is the alternating group $A_n$,
$n\geq 5$, and $C$ is the
set of elements which are the products of two commuting
transpositions.
\item[{\rm (d)}]
$G=SL(2,F)$, where $F$ is a finite field,
and $C$ is the set of conjugates of the
elements $A,B$ in (\ref{defab}).
\end{itemize}
\end{theo}
{\em Proof.} Lemma 2 implies (a) and (b), lemma 3 implies (d).
To prove (c), we remark that $(12)(n-1,n)$ and $(23)(n-1,n)$ generate
a subgroup $D_3$.
$\Box$

Thus we are led to the following problem:
prove or disprove that for all pairs $(G,C)$,
$[K,G,C]$ is either constant or not an invariant
of finite type.

\noindent{\bf Acknowledgement}
\newline
This research is supported by Fonds National Suisse de la
Recherche Scientifique.

\end{document}